\documentclass {amsproc}
\usepackage{amsmath,amsfonts,amsthm,amssymb,hyperref}
\usepackage{amscd}

\usepackage{array}
\usepackage{xy}
\usepackage{epsfig}
\usepackage{graphicx}
\usepackage{color}

\newtheorem{theorem}{Theorem}[section]

\newtheorem {corollary}  [theorem]    {Corollary}

\newtheorem {open}{Open problem}

\theoremstyle{definition}

\newtheorem{example}[theorem]{Example}

\theoremstyle{remark}
\newtheorem{remark}[theorem]{Remark}

\numberwithin{equation}{section}

\def\Z{{\mathbb Z}}

\def\phi{\varphi}

\def\phi{\varphi}

\def\Mg{{\mathfrak g}}

\def\a{\alpha}

\begin{document}

\title[Dimensions of root spaces of hyperbolic Kac--Moody algebras]{Dimensions of Imaginary Root Spaces of Hyperbolic Kac--Moody Algebras}

\author{Lisa Carbone}
\address{Department of Mathematics \\ Rutgers University \\ Piscataway, NJ 08854-8019.}
\email{ carbonel@math.rutgers.edu}
\author{Walter Freyn}
\address{Tu Darmstadt \\ Schlossgartenstr. 7\\ 64289 Darmstadt, Germany. }
\email{freyn@mathematik.tu-darmstadt.de}
\author{Kyu-Hwan Lee}
\address{Department of Mathematics \\ University of Connecticut \\  Storrs, CT 06269-3009.}
\email{khlee@math.uconn.edu}

\begin{abstract}
We discuss the known results and methods for determining root multiplicities for hyperbolic Kac--Moody algebras.
\end{abstract}

\maketitle

\section{Introduction}

Let $\mathfrak{g}$ be a Kac--Moody Lie algebra  over the field $\mathbb C$ of complex numbers, and let $G$ be a Kac--Moody group associated to $\mathfrak{g}$.

 Kac--Moody algebras were introduced in the 1960's independently by Kac,  Moody and Kantor as a generalization of finite dimensional simple Lie algebras. In general, Kac--Moody algebras are infinite dimensional. 
  While they share many structure properties with their finite dimensional predecessors, there appear important new phenomena. Some of the most mysterious ones are related to the \emph{imaginary roots}.  

 There are certain subclasses of special interest. If $\mathfrak{g}$ is of finite type, then $\mathfrak{g}$ is a finite dimensional simple Lie algebra, and $G$ is a simple Lie group. The best-understood infinite dimensional subclass is the one of affine Kac-Moody algebras.  Affine Kac--Moody groups and algebras  give rise to a rich mathematical theory; they are  relevant to number theory and modular forms, lattices and conformal fields theories.

The theory of hyperbolic Kac--Moody groups and algebras naturally generalizes the theory of affine Kac--Moody groups and algebras. Recently hyperbolic and Lorentzian Kac--Moody groups and algebras have been discovered as symmetries in high-energy physics, and they have been shown to serve as duality symmetries  of  a proposed  theory, known as {\it M-theory}, which  unifies all superstring theories (\cite{DHN1,DHN2,HPS,Ju,Kl1,Nic,RW1,RW2,SW1,SW2,We}).

Let $\mathfrak{g}$ be a  Kac--Moody algebra with Cartan subalgebra $\mathfrak{h}$ and root data $\Delta$. The root space $\mathfrak{g}_{\alpha}$ is defined by
$$\mathfrak{g}_{\alpha}\ =\ \{x\in \mathfrak{g}\mid[h,x]=\alpha(h)x,\ h\in \mathfrak{h}\}.$$ Then we have the root space decomposition
$$\mathfrak{g}=\bigoplus_{\alpha\in\Delta^+}\mathfrak{g}_{\alpha}\ \oplus\ \mathfrak{h}\ \oplus\ \bigoplus_{\alpha\in\Delta^-}\mathfrak{g}_{\alpha},$$  
which is a decomposition of $\mathfrak g$ into finite dimensional subspaces, where $\Delta^+$ (resp. $\Delta^-$) is the set of positive (resp. negative) roots. 
The dimension of the root space $\mathfrak{g}_{\alpha}$ is called the {\it multiplicity} of $\alpha$.

Root multiplicities are fundamental data to understand the structure of a Kac--Moody algebra $\mathfrak g$. However, the status of our knowledge shows a dichotomy according to types of $\mathfrak g$. 

Let us recall the classification of types of Kac--Moody algebras $\mathfrak g$. Let $A$ denote the generalized Cartan matrix of $\mathfrak g$. Then the types can be classified as follows.

{\bf Finite type:} $A$ is positive-definite. In this case  $\det(A) > 0$ and $A$ is the Cartan matrix of a finite
dimensional semisimple Lie algebra.

{\bf Affine type:}  $A$ is positive-semidefinite, but not positive-definite. In this case we have $\det(A) = 0$.

{\bf Indefinite type:} $A$ is neither of finite nor affine type. 

The simplest indefinite type is particularly important:

 {\bf Hyperbolic type:} $A$ is neither of finite nor affine type, but every proper, indecomposable principal submatrix is either of finite or affine type. In this case we have $\det(A) < 0$. 

{\bf Lorentzian type:} $\det(A) < 0$ and $A$ has exactly one negative eigenvalue.

The class of Lorentzian generalized Cartan matrices includes but is larger than the class of
hyperbolic generalized Cartan matrices.

Assume that $A$ is of hyperbolic type. Then $A$ is said to be of {\em compact hyperbolic type} if every proper, indecomposable principal submatrix is of finite type.
If $A$ contains an affine submatrix, then it is said to be of {\em noncompact hyperbolic type}.

Now we recall that the associated Weyl group $W$ acts on the set $\Delta$ of all roots, preserving root multiplicities. 
If $\alpha$ is a real root, $\alpha$ has an expression $\alpha=w\alpha_i$ for $w \in W$ where $\alpha_i$ is a simple root. It follows that $\dim(\mathfrak{g}_{\alpha})=1$. Since all roots in finite dimensional Lie algebras are real, all root spaces in finite dimensional Lie algebras are $1$ dimensional. 
 Let $\mathfrak{g}$ be an nontwisted affine Kac--Moody algebra of rank $\ell+1$. Then the multiplicity of every imaginary root of  $\mathfrak{g}$ is $\ell$ (\cite[Corollary 7.4]{K}). There is a similar formula for twisted affine Kac--Moody algebras as well.

For hyperbolic and more general indefinite Kac--Moody algebras the situation is vastly different, due to exponential growth of the imaginary root spaces.
Our knowledge of the dimensions of imaginary root spaces is far from being complete.

\begin{open}
\label{open:find_root_multiplicities} 
{Find an effective closed form formula for the  dimensions of the imaginary root spaces for hyperbolic and other indefinite Kac--Moody algebras.}
\end{open}

A weaker, nevertheles still usefull form of problem~\ref{open:find_root_multiplicities} is the following open question:

\begin{open}
Prove effective upper or lower bounds for the dimensions of the imaginary root spaces for hyperbolic and other indefinite Kac--Moody algebras.
\end{open}

A first step towards this aim might be:

\begin{open}
Conjecture effective upper or lower bounds for the dimensions of the imaginary root spaces for hyperbolic and other indefinite Kac--Moody algebras.
\end{open}

The purpose of this article is to review known results and methods for hyperbolic cases.

\bigskip

The following  gives us a coarse upper bound for each imaginary root space $\Mg_{\a}$, $\a\in\Delta_+^{im}$.

\begin{theorem} [\cite{Ca}] Let $\Mg$ be a Kac--Moody algebra. If $\a\in\Delta_+^{im}$ and $\dim(\mathfrak{g}_{\alpha})>1$, then for $x\in \mathfrak{g}_{\alpha}$, we have
$$x=\sum k_{i_1i_2\dots i_r} [e_{i_1},[e_{i_2},[\dots ,[e_{i_{r-1}},e_{i_r}] ]]$$
where $k_{i_1i_2\dots i_r}\in\mathbb{C}$, $\a=\alpha_{i_1}+\alpha_{i_2} +\dots +\alpha_{i_{r}}$ for simple roots
$\alpha_{i_1}$, $\alpha_{i_2}$, $\dots$, $\alpha_{i_{r}}$,  not necessarily distinct, and the sum is taken over the $\ell^r$ orderings of 
$\{\a_{i_1},\a_{i_2},\dots ,\a_{i_r}\}$, with all but $\dim(\Mg_{\a})$ of the $k_{i_1i_2\dots i_r}$ equal to 0.
\end{theorem}

\begin{corollary} [\cite{K}] Let $\Mg$ be  any Kac--Moody algebra. Let $\a\in\Delta^{im}$. Then
$$\dim(\mathfrak{g}_{\alpha})\leq \ell^{|height(\alpha)|}<\infty.$$
\end{corollary}

\section{Berman-Moody formula and Peterson formula}

 The first formulas for root multiplicities of Kac--Moody algebras are a closed form formula by Berman and Moody (\cite{BM}) and a recursive formula by Peterson (\cite{P}). Both formulas are based on the {\it denominator identity} for a Kac--Moody algebra $\mathfrak{g}$:
$$\prod_{\alpha\in\Delta^+} (1-e(-\alpha))^{\rm{mult}(\alpha)}\quad =\quad 
\sum_{w\in W}(-1)^{l(w)}e(w(\rho)-\rho),$$
where $\rm{mult}(\alpha)=\dim(\mathfrak{g}_{\alpha})$, $e(\alpha)$ is a formal exponential, $l(w)$ is the length function on $W$ and $\rho \in\mathfrak{h}^{\ast}$ satisfies $\langle\rho,\alpha_i^{\vee}\rangle=1$, $i=1,\dots \ell$.
The denominator identity is derived from the {\it Weyl-Kac character formula}  for the trivial one dimensional module, which has character equal to 1, and relates the orbits of the Weyl group $W$ and the root multiplicities of the Kac--Moody algebra.

 The simple roots determine the Weyl group and hence the right hand side of the denominator identity. The root multiplicities determine the left hand side of the denominator identity. Therefore if we know the simple roots, the denominator identity can in principle be used to give an alternating sum formula for the root multiplicities. This is not effective in general, but as a systematic application, Peterson (\cite{P}) used the denominator identity to give a recursive formula for the root multiplicities of indefinite algebras. Peterson's paper has never been published. However, an outline of the proof can be found in Exercises 11.11 and 11.12 in \cite{K}. One can  find a complete proof by Kang, Kwon and Oh \cite{KKO} in a more general setting of graded Lie superalgebras. There is an improved version of Peterson's formula by Moody and Patera \cite{MP}. 

Implementations of Peterson's formula have been used by many authors. For example, Peterson's formula was used  to produce the root multiplicity table in \cite{K} for the hyperbolic Kac--Moody algebra $\mathfrak F$ with Cartan matrix  $\begin{pmatrix} 2 & -2 & 0 \\ -2 & 2 & -1 \\ 0 & -1 & 2 \end{pmatrix}$ of type $HA_1^{(1)}$. Borcherds \cite{Bo1} used Peterson's formula to compute root multiplicities of the Lorentzian Kac--Moody algebra attached to the even unimodular lattice $II_{25,1}$. B\"arwald, Gebert and Nicolai \cite{BGN} used the formula for $E_{10}$, and Kleinschmidt \cite{Kl2} used it to produce  extensive tables for root multiplicities of $E_{10}$ and $E_{11}$.

Berman and Moody's formula was also derived from the denominator identity by taking logarithmic derivatives and using M{\"o}bius inversion. The resulting formula is in a closed form. To describe this formula, we introduce the following notation. We write $\lambda | \alpha$ if $\alpha=r \lambda$ for some positive integer $r$ and denote $1/r$ by $\lambda / \alpha$. For each $w \in W$, let $s(w)$ be the sum of the positive roots which are mapped into $\Delta^-$ by $w^{-1}$, and  for each $s(w)$, set 
$\varepsilon(s(w))= (-1)^{l(w)+1}$. We enumerate $s(w)$, $1 \neq w \in W$ as $s_1, s_2, \dots$ in an order of increasing height. For $\lambda \in Q:= \oplus \mathbb Z_{\ge 0} \alpha_i$, the set $S(\lambda)$ is defined to be $\{ (n)=(n_1, n_2, \dots ) | \sum n_i s_i =\lambda \}$. Then the Berman-Moody formula is given by
\begin{equation} \label{BM} \mathrm{mult} (\alpha)= \sum_{\lambda | \alpha} \mu \left (\frac {\alpha}{\lambda} \right ) \frac {\lambda}{\alpha} \sum_{(n) \in S(\lambda)} \prod \varepsilon (s_i)^{n_i} \ \frac{( (\sum n_i) -1)!}{\prod (n_i!)} ,\end{equation}
where $\mu$ is the M{\"o}bius function.

\begin{example}
Let $A=\begin{pmatrix} 2 & -3 \\ -3 &2 \end{pmatrix}$, and consider $\alpha= 4 \alpha_1 + 5 \alpha_2$. Let $w_1$ and $w_2$ be the simple reflections. We set \[ \begin{array}{ll} s_1=s(w_1)=\alpha_1, & s_2=s(w_2)=\alpha_2, \\ s_3=s(w_1w_2)=4\alpha_1 +\alpha_2, & s_4=s(w_2w_1)=\alpha_1 +4 \alpha_2 .\end{array} \] Then we obtain \[ S(\alpha) = \{(4,5,0,0),(3,1,0,1), (0,4,1,0) \}, \]
and get \[ \mathrm{mult}(\alpha) = \frac {8!}{4!5!} -\frac{4!}{3!} -\frac{4!}{4!} = 14 -4-1 =9.\]
\end{example}

Berman and Moody's formula \eqref{BM} can be considered as a generalization of Witt formula \cite{Se} for free Lie algebras and  was later generalized by Kang (\cite{Ka1, Ka2}). 

Both Peterson's formula and Berman-Moody's formula enable us to calculate the multiplicity of a given root (of a reasonable height), but they do not provide any insight into properties of root multiplicities.

\bigskip

\section{Works of Feingold and Frenkel, and Kac, Moody and Wakimoto}

Pioneering works on root multiplicities of hyperbolic Kac--Moody algebras began with the paper by Feingold and Frenkel \cite{FF}, where the hyperbolic Kac--Moody algebra $\mathfrak F$ of type $HA_1^{(1)}$ was considered. Using the same method,  Kac, Moody and Wakimoto \cite{KMW} calculated some root multiplicities for $HE_8^{(1)}(=E_{10})$.

Their method is based on the fact that any symmetrizable Kac- Moody Lie algebra can be realized as the minimal graded Lie algebra $\mathfrak g = \bigoplus_{n\in \mathbb Z} \mathfrak g_n$ with the local part $V\oplus \mathfrak g_0 \oplus V^*$, where $\mathfrak g_0$ is a smaller (typically,  finite-dimensional or affine) Kac--Moody Lie algebra, $V$ is an integrable irreducible highest weight representation of $\mathfrak g_0$, and $V^*$ is the contragredient of $V$. More precisely, we consider the graded Lie algebra $G= \bigoplus_{n\in \mathbb Z} G_n$, where $G_0=\mathfrak g_0$ and  $G_{\pm} = \bigoplus_{n \ge 1} G_{\pm n}$ be the free Lie algebra generated by $G_1=V^*$ (respectively, $G_{-1}=V$). We let $I = \bigoplus_{n\in \mathbb Z} I_n$ be the maximal graded ideal of $G$ intersecting the local part $V\oplus \mathfrak g_0 \oplus V^*$ trivially. Then we obtain $\mathfrak g=G/I$. 

Since $G_-$ is a free Lie algebra, its homogeneous dimensions can be computed using  the generalized Witt formula \cite{Ka1} and one can use  the representation theory of $\mathfrak g_0$ to determine the structure of $I_n$ and get root multiplicities for $\mathfrak g_n$. However, if $n$ is big, it is not easy to handle $I_n$ explicitly.

\begin{example}
We again consider $A=\begin{pmatrix} 2 & -3 \\ -3 &2 \end{pmatrix}$ and $\alpha= -4 \alpha_1 - 5 \alpha_2$. We keep the notations in the previous paragraph. We take the subalgebra $\mathfrak g_0$ generated by the Chevalley generators $e_2$, $f_2$ and the Cartan subalgebra $H=\mathbb C h_1 + \mathbb C h_2$. Then we have $\mathfrak g_0 \cong sl_2(\mathbb C) \oplus \mathbb C h_1$. Consider $f_1$ as a highest weight vector for $\mathfrak g_0$ under adjoint action.
Since $[h_2, f_1]=3f_1$, the element $f_1$ generates the $4$-dimensional irreducible representation $V$ with a basis $\{ f_1, [f_1f_2], [[f_1 f_2]f_2], [[[f_1f_2]f_2]f_2]\}$. From the Serre relation, we see that $I_{-i}=0$ for $0 \le i \le 3$. Since $(-4\alpha_1-\alpha_2)(h_2)=10$, the graded piece $I_{-4}$ is isomorphic to the $11$-dimensional irreducible $\mathfrak g_0$-module generated by $[f_1[f_1[f_1[f_1f_2]]]]$. From the representation theory of $sl_2(\mathbb C)$, we know that each weight space of $I_{-4}$ is $1$-dimensional. Using the generalized Witt's formula for free Lie algebras (cf. \cite{Ka1}), we compute the dimension of the subspace of $G_{-4}$ with weight $-4 \alpha_1 - 5 \alpha_2$ and obtain $10$. Since the dimension of the subspace of $I_{-4}$ with weight $-4 \alpha_1 - 5 \alpha_2$ is $1$ as already observed, the root multiplicity of $-4 \alpha_1 - 5 \alpha_2$ for $\mathfrak g$ is $10-1=9$.

\end{example}

Using this inductive  method, in the papers \cite{FF, KMW}, the authors calculated root multiplicities up to level $2$; that is, multiplicities for the roots belonging to $\mathfrak g_1$ and $\mathfrak g_2$. The results involve various partition functions, and will be presented in what follows.

First, we recall that the rank 3 Kac--Moody algebra $\mathfrak F$ of type $HA^{(1)}_1$  has Weyl group isomorphic to $PGL_2({\Z})$ and that positive imaginary roots of  $\mathfrak F$ are identified with $2 \times 2$ semi-positive definite symmetric matrices. The algebra $\mathfrak F$ has noncompact hyperbolic type  and was first studied by Feingold and Frenkel \cite{FF}. They choose a subalgebra $\mathfrak F_0$ of affine type $A^{(1)}_1$ as the smaller algebra in the inductive method. Then the level 1 root multiplicities of $\mathfrak F$ are obtained from the weight multiplicities of the basic representation of the affine algebra $\mathfrak F_0$. The character of the basic representation was calculated by Feingold and Lepowsky \cite{FL}. Using the character, it was shown in \cite{FF} that for each level 1 root $\alpha$, we have
$$\dim({\mathfrak F}^{\alpha})= p \left(1-\frac{(\alpha |\alpha)}{2}\right),$$
where $p(n)$ is the classical partition function.

 They also obtained the generating function of the multiplicities of the level $2$ roots: 
\begin{eqnarray} & & \sum_{n=0}^\infty M(n-1) q^n \nonumber \\ &=& \frac {q^{-3}}2 \left(\sum_{n=0}^\infty p(n)q^n \right) \prod_{j=1}^\infty (1-q^{4j-2}) \nonumber \\ & & \qquad \qquad \times \left(\prod_{j=1}^\infty (1+q^{2j-1})-\prod_{j=1}^\infty (1-q^{2j-1})-2q \right) \nonumber \\
&=& (1 -q^{20} + q^{22}-q^{24}+q^{26}-2 q^{28} + \cdots) \sum_{n=0}^\infty p(n) q^n , \label{eqn-ff}
\end{eqnarray}
where $M(2m) = \mathrm{mult}  \begin{pmatrix} m & 0 \\ 0 &2 \end{pmatrix}$ and $M(2m-1)= \mathrm{mult}  \begin{pmatrix} m &1 \\ 1 & 2 \end{pmatrix}$ with identification of  positive imaginary roots with $2 \times 2$ semi-positive definite symmetric matrices. In particular, this supports Frenkel's conjecture in Section \ref{sec-F}.

 Next, let $p^{(\ell)}(n)$ be the number of partitions of $n$ into parts of $\ell$ colors. We write $\phi(q)=\prod_{n=1}^\infty (1-q^n)$. Then the generating function of $p^{(\ell)}(n)$ is given by 
\begin{equation} \label{pl} \sum_{n=0}^{\infty} p^{(\ell)}(n) q^n=\dfrac{1}{\phi(q)^{\ell}}. \end{equation}
Kac, Moody and Wakimoto \cite{KMW} used the inductive method to calculate root multiplicities of the hyperbolic Kac--Moody algebra $\mathfrak g$ of type $E_{10}$ for the roots of level $\le 2$. They showed
\begin{equation} \label{eqn-KMW} \dim({\mathfrak {g}}^{\alpha})=\begin{cases} p^{(8)}\left(1-\frac{(\alpha |\alpha)}{2}\right) & \text{ if } \alpha \text{ is of level 0 or 1},\\ \xi\left(3-\frac{(\alpha |\alpha)}{2}\right) & \text{ if } \alpha \text{ is of level 2},
\end{cases} \end{equation} where the function $\xi(n)$ is given by
\[ \sum_{n=0}^\infty \xi(n) q^n = \dfrac{1}{\phi(q)^{8}} \left [ 1 - \dfrac {\phi(q^2)}{\phi(q^4)} \right ]. \]

\bigskip 


\section{Kang's generalizations}

The inductive method used in \cite{FF, KMW} was systematically developed and generalized by Kang \cite{Ka1} for arbitrary Kac--Moody algebras and for higher levels, and has been adopted in many works on roots multiplicities of indefinite Kac--Moody algebras. In his construction, Kang adopted homological techniques and Kostant's formula (\cite{GL}) to devise a method that works for higher levels. For example, Kang applied his method to compute roots multiplicities of the algebra $\mathfrak F$ of type $HA_1^{(1)}$ up to level $5$ (\cite{Ka3, Ka4}).

Kang also took a different approach and generalized Berman and Moody's formula in \cite{Ka2}. The resulting formula does not require the generalized Witt formula for free Lie algebras. As in the inductive method, he chooses a smaller subalgebra $\mathfrak g_0$ of $\mathfrak g$ and considers the decomposition 
\[ \mathfrak g= \mathfrak g_- \oplus \mathfrak g_0 \oplus \mathfrak g_+ ,  \] where $\mathfrak g_- = \sum_{n<0} \mathfrak g_n$ and $\mathfrak g_+ = \sum_{n>0} \mathfrak g_n$. He then uses
the chain complex of $\mathfrak g_-$ and its homology modules to apply the Euler-Poincar{\'e} principle and obtains an identity:
\[ \sum_{k=0}^\infty (-1)^k \mathrm{ch} \ \Lambda^k(\mathfrak g_-) = \sum_{k=0}^\infty (-1)^k \mathrm{ch}\ H_k(\mathfrak g_-), \] where $\mathrm{ch}\ V$ is the character of a $\mathfrak g_-$-module $V$. This identity can be considered as a generalized denominator identity for the algebra $\mathfrak g$ with respect to $\mathfrak g_0$.  From this identity, Kang derives a recursive formula (\cite[Theorem 3.3]{Ka2}) and a closed form formula (\cite[Theorem 4.1]{Ka2})  for all root multiplicities, though an actual computation, in general,  requires a substantial amount of information on representations of the Lie algebra $\mathfrak g_0$. The closed form formula can be considered as a generalization of Berman and Moody's formula. By choosing various $\mathfrak g_0$, one can have more flexibility in computations and obtain interesting identities. 

In order to understand the representations of $\mathfrak g_0$ when $\mathfrak g_0$ is of affine type, Kang adopted the path realization of crystal basis (\cite{KMN1,KMN2}) for affine Kac--Moody algebras in  \cite{Ka3,Ka4,KM1}. This idea was followed by Klima and Misra \cite{KMi} for the indefinite Kac--Moody algebras of symplectic type.

 There are many partial results for root multiplicities of hyperbolic and Lorentzian Kac-Mooy algebras by applying Kang's methods. 
 As mentioned earlier,  Kang  used his inductive method to give root multiplicity formulas for the roots of $\mathfrak F $ of level $\leq$ 5 \cite{Ka3, Ka4}.  In an unpublished work, Kac also discovered a level 3 root multiplicity formula for $\mathfrak F$. The inductive method can be applied for higher level roots, but it would be a daunting task to derive any concrete formula as the complexity of computations grows fast in this method. Furthermore, even for levels $3$ and $4$,  the generalized Witt formula requires long computations of partitions for roots of large heights. Thus the following problem is wide open.

 \begin{open}{ Find a root multiplicity function for roots of $\mathfrak F$ of all levels. }
\end{open}


 In addition to the works on $\mathfrak F$, Kang's methods were used in the papers by  Benkart, Kang and Misra \cite{BKM1, BKM2, BKM3},  Kang and Melville \cite{KM1,KM2}, Hontz and Misra \cite{HM1, HM2}, Klima and Misra \cite{KMi}.
 In order to discuss these results on root multiplicities,  we introduce a notation for generalized Cartan matrices, following \cite{BKM2, BKM3}: we define

$$IX(a,b)=\left(\begin{matrix}
{2} & {-b} & {0} & {\cdots} & 0  \\
{-a} & {} & {} & {} & {} \\
{0} & {} &  & & {}  \\
{\vdots} & {} &  C(X) & & {} \\
0 & {} & {}  & {} & {} \\
\end{matrix}
\right)$$
where $C(X)$ is the Cartan matrix of a Kac--Moody algebra of type $X$. For almost all positive integer values of $a$ and $b$, $IX(a,b)$ is a generalized Cartan matrix of indefinite type.

 $\circ$ $IA_1(a,b)=\left( \begin{matrix}
2 & -b\\
-a & 2
\end{matrix}
\right)$: This is a rank 2 generalized Cartan matrix which is hyperbolic if ${ab > 4}$. Its connection to Hilbert modular forms was observed in an earlier paper by Lepowsky and Moody \cite{LM}.  The root multiplicities of this algebra were studied in \cite{BKM1, KM2}. In particular, Kang and Melville \cite{KM2} made various observations on the relationships between root length and root multiplicity for this algebra.

$\circ$ $IA_n(a,b)$ and $IX_n(a,1)$, $X=B,C,D$: Benkart, Kang and Misra \cite{BKM2, BKM3} studied the associated indefinite Kac--Moody algebras and found formulas for root multiplicities for roots of degree $\le 2a+1$ for type $A$ and of degree $\le 2a$ for other types. The multiplicity formulas are obtained by exploiting representation theory of finite dimensional simple Lie algebras and involve the Littlewood-Richardson coefficients and Kostka numbers.  
 
$\circ$ $IA_n^{(1)}(1,1)=HA_n^{(1)}$: Kang and Melville \cite{KM1} considered the root multiplicities of the hyperbolic Kac--Moody algebra of this type. ( Notice that if $n=1$ then we obtain the algebra $\mathfrak F$.) They gave a formula for all root multiplicities  using the path realization of affine crystals \cite{KMN1, KMN2}, but the formula depends on how to compute the path realization. Hontz and Misra \cite{HM1} also calculated some root multiplicities  for this type.

$\circ$ $ID_4^{(3)}(1,1)=HD_4^{(3)}$ and $IG_2^{(1)}(1,1)=HG_2^{(1)}$: Hontz and Misra \cite{HM2} considered these types. They chose a subalgebra $\mathfrak{gl}(4, \mathbb C)$ to use the inductive method and the combinatorics of the representations of $\mathfrak{gl}(4, \mathbb C)$. They obtained root multiplicities for roots of degree $\le 8$.

$\circ$ $IC_n^{(1)}(1,1)=HC_n^{(1,1)}$: Klima and Misra \cite{KMi} used Kang's formula and the path realization of the crystals to obtain some root multiplicities  and showed that the multiplicities of these roots are polynomials in $n$.

The above results exploited various constructions and methods in the representation theory of finite-dimensional and affine Kac--Moody algebras in order to apply the inductive method.
 Root multiplicities of other indefinite algebras were also obtained by \cite{SUL} and \cite{SLU}. Kang's methods have been further developed to include the cases of generalized Kac-Moody algebras \cite{Ka5} and more general graded Lie (super)algebras \cite{Ka6}.


\bigskip
\section{Frenkel's conjecture} \label{sec-F}

Despite all the results in the previous sections, we still do not have any unified, efficient approach to computing all root multiplicities or  explicit bounds. Essentially, the methods give answers for root multiplicities one at a time, with no general formulas or effective bounds on multiplicities. Furthermore, as far as we know, there is no conjecture on root multiplicities for general indefinite Kac--Moody algebras.

For hyperbolic Kac--Moody algebras, in the setting of the `no-ghost' theorem from string theory, I.  Frenkel \cite{F} proposed a bound on the root multiplicities of hyperbolic Kac--Moody algebras.

\medskip

{\bf Frenkel's conjecture:} {\it Let ${\mathfrak {g}}$ be a symmetric hyperbolic Kac--Moody algebra associated to a hyperbolic lattice of dimension $d$ and equipped with invariant form $(\cdot\mid\cdot)$ normalized to equal 2 on simple roots. Then we have:
$$\dim({\mathfrak {g}}^{\alpha})\ \leq\  p^{(d-2)} \left(1-\frac{(\alpha |\alpha)}{2}\right),$$
where the function $p^{(\ell)}(n)$  is defined in \eqref{pl}. }

\begin{remark}
A caveat is that the dimension of a hyperbolic lattice and the rank of the corresponding Kac--Moody algebra may be different. For example, the dimension of the even unimodular lattice $II_{25,1}$ is 26, but the rank of the corresponding Kac--Moody algebra is infinite. This fact is related to Conway's work on $II_{25,1}$ \cite{Con}.
\end{remark}

Frenkel's conjecture is known to be true in a number of important cases. In particular, there is a distinguished symmetric Lorentzian Kac--Moody algebra $\mathfrak{L}$ of rank $26$. The Dynkin diagram of $\mathfrak{L}$ corresponds to the Dynkin diagram of the $26$ dimensional even unimodular Lorentzian lattice $II_{25,1}$. Conway \cite{Con} showed that the Dynkin diagram is given by the Leech lattice, i.e. the unique $24$-dimensional positive-definite even unimodular lattice with no roots.
In \cite{F}, Frenkel used the no-ghost theorem to show that 
$$\dim({\mathfrak {L}}^{\alpha})\ \leq\  p^{(24)} \left(1-\frac{(\alpha |\alpha)}{2}\right).$$ Frenkel's method gives the same upper bound for the multiplicities of any symmetric Kac--Moody algebra associated to a hyperbolic lattice of dimension $26$,  the {\it critical dimension}.

Borcherds computed the root multiplicities of $\mathfrak{L}$ using Peterson's recursive formula (\cite{Bo1}), and observed that  `most' root vectors of small norm have root multiplicities {\it equal} to Frenkel's upper bound. He then constructed a {\em generalized} Kac--Moody algebra  whose root multiplicities are exactly the same as the Frenkel's bound. (See Section \ref{sec-B}.) 

The conjecture is still open for the rank $3$ hyperbolic Kac--Moody algebra $\mathfrak F$.

\begin{open} \cite[Exercise 13.37]{K}
Prove Frenkel's cojecture for the rank $3$ hyperbolic algebra $\mathfrak F$.
\end{open}

 Klima and Misra \cite{KMi} showed that Frenkel's bounds do not hold for  indefinite Kac--Moody algebras of symplectic type.  However, this case is not of type ADE and hence is not included in Frenkel's conjecture.


\section{Extended and overextended Dynkin diagrams} \label{extended}

Let $\Delta$ be a finite indecomposable root system,  that is, the Dynkin diagram of a root system of a finite dimensional Lie algebra.  Let $\Pi=\{\alpha_1,\dots ,\alpha_{\ell}\}$ be the simple roots of $\Delta$. For $\Delta$ indecomposable, there is a unique root $\theta$ called the {\it highest  root} that is a linear combination of the simple roots with positive integer coefficients. The highest root $\theta$  satisfies $(\theta,\alpha)\geq 0$ for every simple root $\alpha$ and $(\theta,\beta)> 0$ for some simple root $\beta$, where $(\cdot,\cdot)$ is the positive definite symmetric bilinear form corresponding to $\Delta$.

Let $\Pi'=\Pi\cup\{\alpha_0\}$. Then $\Pi'$ is called the {\it extended system of simple roots} corresponding to $\Pi$.   Let $\delta=\alpha_0+\theta$, and $\delta$ is the minimal {\em null root} such that $(\delta,\delta)=0$.
The Dynkin diagram of $\Pi'$ is called the {\it extended Dynkin diagram} or {\em untwisted affine Dynkin diagram} corresponding to $\Delta$. An extended Dynkin diagram has a vertex labeled 0 corresponding to the root $\alpha_0$. 

A generalized Cartan matrix $A$ is called {\it Lorentzian} if $\det(A)\neq 0$ and $A$ has  exactly one negative eigenvalue. A  {\it Lorentzian Dynkin diagram} is the Dynkin diagram of a  Lorentzian generalized Cartan matrix.  A {\it Lorentzian extension} $\mathcal{D}$ of an untwisted affine Dynkin diagram $\mathcal{D}_0$ is a Dynkin diagram obtained by adding one vertex, labeled $-1$, to $\mathcal{D}_0$ and connecting the vertex $-1$ to the vertex of $\mathcal{D}_0$ labeled 0 with a single edge. A common notation for the Lorentzian extension of some finite simple Dynkin diagram $X$ is $X^{++}$. It is sometimes called the {\em over-extended} diagram. The addition of a third additional vertex then is denoted by $X^{+++}$.

 Every Lorentzian extension of an untwisted affine Dynkin diagram is a Lorentzian Dynkin diagram, in fact a hyperbolic Dynkin diagram. 

\medskip
\begin{example}[$E_{10}$ and $E_{11}$] Let $\Delta$ be the Dynkin diagram for $E_8$. We label the first vertex of the `long tail' by 1. Adding a vertex labeled 0 and connecting vertices 0 and 1 by a single edge yields the extended Dynkin diagram $\Delta'$ which corresponds to the affine Kac--Moody algebra $E_9=E_8^{(1)}$. Adding a further vertex labeled $-1$ and connecting vertices $-1$ and 0 by a single edge yields the overextended Dynkin diagram  which corresponds to the hyperbolic Kac--Moody algebra $E_{10}$. Adding still another vertex yields $E_{11}$. Hence we have 
$$E_8^+=E_9,\qquad E_8^{++}=E_{10},\qquad E_{8}^{+++}=E_{11}\,.$$

\end{example}

\section{Root multiplicities for $E_{10}$ and $E_{11}$}

The Lie algebras of types $E_{10}$ and $E_{11}$ are of special importance due to various conjectures describing their appearance in string theory and $M$-theory. $E_{10}$ also is distinguished due to results of Viswanath, showing, that any simply-laced hyperbolic Kac--Moody algebra embeds into $E_{10}$ (\cite{Vis}).
 
Bauer and Bernard \cite{BB} found the root multiplicities of $E_{10}$ and $\mathfrak F$ up to level $3$ using the inductive method. However, their computations are written in terms of {\em conformal blocks} rather than {\em partition functions}.

Even though Frenkel's conjecture provides an important guideline to root multiplicities of hyperbolic Kac--Moody algebras, there are some known counterexamples.

For example,  Kac, Moody and Wakimoto \cite{KMW} showed that the conjecture fails for $E_{10}$. More precisely, using the notations in \eqref{eqn-KMW}, we have  $$\xi(6)=p^{(8)}(4)+1>p^{(8)}(4).$$ This implies that $E_{10}$ does not satisfy Frenkel's conjecture at level $2$.

Extensive calculations of root multiplicities of over-extended, simply-laced hyperbolic Kac--Moody algebras were done by A.~Kleinschmidt in his thesis \cite{Kl1}. His calculations suggest that for the Lorentzian extensions of $A_n$ and $D_n$, the bounds given by Frenkel's conjecture hold. Nevertheless Frenkel's bounds  fail for $E_8^{++}=E_{10}$. If we add one more vertex to the Lorentzian extensions, the behavior of root multiplicities seems to be again in agreement with Frenkel's conjecture. More unpublished calculations by A.~Kleinschmidt and H.~Nicolai confirm these findings, that is, the  failure of  Frenkel's conjecture for $E_{10}$ at higher levels (\cite{KN}), and 
validity of Frenkel's conjecture for $E_{11}$ (\cite{KN}). Large tables of root multiplicities of $E_{10}$ and $E_{11}$ can be found in \cite{Kl2}.


\section{Borcherds's constructions} \label{sec-B}

As mentioned earlier, Frenkel showed that the root multiplicities of the Lorentzian Kac--Moody algebra $\mathfrak{L}$ has an upper bound  $p^{(24)}(1-(\alpha|\alpha)/2)$. We note that the function $p^{(24)}(n)$ is related to the modular discriminant $\Delta(z)$ in the theory of modular forms, where $\Delta(z)=q\prod_{n=1}^\infty (1-q^n)^{24}=q\phi(q)^{24}$, $q=e^{2\pi i z}$ for $z \in \mathcal H$ the upper half plane. 

One can ask: Is it possible to construct a Lie algebra whose root multiplicities are exactly $p^{(24)}(1-(\alpha|\alpha)/2)$? 

Such a Lie algebra would show `modular behavior'. Indeed, in \cite{Bo3},  Borcherds constructed the {\em fake Monster Lie algebra} $\mathfrak {M}$ which contains the Kac--Moody algebra $\mathfrak {L}$ and whose root multiplicities are exactly given by $p^{(24)}(1-(\alpha|\alpha)/2)$. 

The fake Monster Lie algebra $\mathfrak M$ was obtained from a {\em lattice vertex algebra} \cite{Bo2}. Let $M$ be a nonsingular even lattice, and $V(M)$ be the vertex algebra associated to $M$. Then we have the Virasoro operators $L_i$ on $V(M)$ for each $i \in \mathbb Z$. We define the {\em physical space} $P^n$ for each $n \in \mathbb Z$ to be the space of vectors $w \in V(M)$ such that
\[  L_0(w) = nw \qquad \text{ and } \qquad L_i(w)=0 \ \text{ for } i >0 .\] 
Then the space $\mathfrak G_M:= P^1/L_{-1} P^0$ is a Lie algebra
 and satisfies the following properties \cite{Bo2}:
\begin{enumerate}
\item Let $\mathfrak g$ be a Kac--Moody Lie algebra with a generalized Cartan matrix $A$ that is indecomposable, simply laced and non-affine. If the lattice $M$ contains the root lattice of $\mathfrak g$ then $\mathfrak g$ can be mapped into $\mathfrak G_M$ so that the kernel is in the center of $\mathfrak g$.

\item Let $d$ be the dimension of $M$, and $\alpha \in M$ be a root such that $(\alpha | \alpha) \le 0$. Then the root multiplicities of $\alpha$ in $\mathfrak G_M$ is equal to \begin{equation} \label{eqn-bound} p^{(d-1)}(1-(\alpha|\alpha)/2) - p^{(d-1)}(-(\alpha|\alpha)/2).\end{equation}
\end{enumerate}

Therefore, when we have a hyperbolic Kac--Moody algebra $\mathfrak g$ with root lattice $M$, the Lie algebra $\mathfrak G_M$ contains $\mathfrak g$ and provides an upper bound \eqref{eqn-bound} for root multiplicities of $\mathfrak g$. Note that this bound is weaker than Frenkel's conjecture.    
The Lie algebra $\mathfrak G_M$ is not a Kac--Moody algebra but a {\em generalized} Kac--Moody algebra or a {\em Borcherds algebra}, since  it has imaginary simple roots.

When $M=II_{25,1}$, the fake Monster $\mathfrak M$ is obtained by taking a quotient of $\mathfrak  G_M$ by the kernel of a bilinear form, and the no-ghost theorem can be utilized to show that the root multiplicities are exactly given by $p^{(24)}(1-(\alpha|\alpha)/2)$. As a result, Borcherds could write down the denominator identity of $\mathfrak M$:
\[ e(\rho) \prod_{\alpha \in \Delta^+} (1-e(-\alpha))^{p^{(24)}(1-(\alpha|\alpha)/2)}=\sum_{\substack {w \in W  \\ n \in \mathbb N}} \det(w) \tau(n) e(w(n\rho)) ,\] where $\tau(n)$ is the Ramanujan tau function defined to be the Fourier coefficients of $\Delta (z)$, i.e.  $\sum_{n \in \mathbb N} \tau(n) q^n = q \prod_{n=1}^\infty (1-q^n)^{24}=\Delta (z)$. Borcherds also showed that the denominator function is itself an automorphic form. 

While describing root multiplicities of an indefinite Kac--Moody algebra is very difficult,  Borcherds constructions produce many examples of generalized Kac--Moody algebras whose root multiplicities are explicitly known. Moreover, Borcherds's examples extend some Kac--Moody algebras to generalized Kac--Moody algebras so that we may obtain automorphic forms from the denominator functions of the generalized Kac--Moody algebras. Pursuing Borcherds's idea, Gritsenko, Nikulin and Niemann constructed generalized Kac--Moody algebras to extend some hyperbolic Kac--Moody algebras \cite{GN1, GN2, Nie}.

For the hyperbolic Kac--Moody algebra $\mathfrak F$, possible connections to Siegel modular forms were noticed by Feingold and Frenkel \cite{FF}. 
Gritsenko and Nikulin's construction indeed shows that  the denominator function of the corresponding generalized Kac--Moody algebra is a Siegel modular form.
More precisely, 
they showed that there exists a generalized Kac--Moody algebra $\mathcal G$ which contains $\mathfrak F$ and whose denominator function is the weight $35$ Siegel cusp form $\Delta_{35}(Z)$, which is called the {\em Igusa modular form}. As a byproduct, they obtained the infinite product expression of $\Delta_{35}(Z)$.
Even though this construction manifests the connection of $\mathfrak F$ to a Siegel modular form, the root spaces of $\mathcal G$ are much bigger than those of $\mathfrak F$, and the construction does not help understand root multiplicities of $\mathfrak F$.

In his Ph.D. thesis \cite{Nie}, P. Niemann constructed a generalized Kac--Moody algebra $\mathcal G_{23}$ which contains $\mathfrak F$. The denominator function of $\mathcal G_{23}$ is closely related to the eta product $\eta(z)\eta(23z)$, where $\eta$ is the Dedekind $\eta$-function. If $q \eta^{-1}(z) \eta^{-1}(23 z)=\sum_{n=0}^\infty p_{\sigma}(n) q^n$, $q=e^{2 \pi i z}$, he showed that
\[ \mathrm{mult}(\mathfrak F, \alpha) \le \begin{cases} p_\sigma ( 1 - \tfrac 1 2 (\alpha , \alpha) ) & \text{ if } \alpha \notin 23 L^* , \\ p_\sigma ( 1 - \tfrac 1 2 (\alpha , \alpha) ) +  p_\sigma ( 1 - \tfrac 1 {46} (\alpha , \alpha) ) & \text{ if } \alpha \in 23 L^* ,
\end{cases}
\]where $L$ is a certain lattice and $L^{*}$ is its dual.
This bound is quite close to Frenkel's conjecture.
Note that we have
\begin{align*} \label{eqn-opop} \sum_{n=0}^\infty p_\sigma(n) q^n &= (1+ q^{23} + 2 q^{46} + \cdots ) \prod_{n=1}^\infty (1-q^n)^{-1} \\ &= (1+ q^{23} + 2 q^{46} + \cdots ) \sum_{n=0}^\infty p(n) q^n .\end{align*}
One can also compare this with the actual multiplicities of level $2$ roots given in \eqref{eqn-ff}.

\bigskip


\section{Asymptotics by the method of Hardy-Ramanujan-Rademacher}

Now that root multiplicities of  $\mathcal G_{23}$ are given by Fourier coefficients of automorphic forms, we can apply analytic tools to get asymptotic formulas for these multiplicities; namely, one can use the method of Hardy-Ramanujan-Rademacher to obtain asymptotic formulas for $p_\sigma(1+n)$. See e.g. \cite{Leh} for the details of this method. 

As for Niemann's bound, we consider $f(z)=\eta(z)^{-1}\eta(23z)^{-1}$, which is a weakly holomorphic modular form of weight $-1$ with respect to
$\Gamma_0(23)$, where \[ \Gamma_0(23) =\left \{ \begin{pmatrix} a & b \\ c & d \end{pmatrix} \in SL_2(\mathbb Z) :  c \equiv 0 (\mathrm{mod}\ 23) \right \}. \]

Kim and Lee obtained the following asymptotics using the method of Hardy-Ramanujan-Rademacher:
\begin{theorem} \cite{KL} \label{niemann}
$$p_{\sigma}(n+1)= \frac {2\pi}{n\sqrt{23}} I_2 \left ( \frac {4\pi\sqrt{n}}{\sqrt{23}} \right ) +O\left(n^{-\frac 12} I_2 \left ( \frac {2\pi\sqrt{n}}{\sqrt{23}} \right )\right),
$$
where $I_2$ is the modified Bessel function of the first kind.
\end{theorem}

This result has an immediate implication on root multiplicities of the hyperbolic Kac--Moody algebra $\mathfrak F$. For example, if $(\alpha, \alpha) = -56$ then the main term of the asymptotic formula gives $4578.99$, while the actual value of the Fourier coefficient is $4576$. The exact value of $\mathrm{mult}(\mathfrak F, \alpha)$ is $4557$. In this way, we can calculate a sharp upper bound for $\mathrm{mult}(\mathfrak F, \alpha)$ even if $|(\alpha, \alpha)|$ is big. 

\begin{example} \label{example}
If $\alpha= 10\alpha_1 + 10 \alpha_2 + 5\alpha_3=(10,10,5)$ then $-\frac 1 2 (\alpha, \alpha) = 25$ and we have 
$\mathrm{mult}(\alpha) = 2434$, and the main term of the asymptotics is $2437.16$. We calculate more cases and make a table:
\begin{equation}
\label{table}
\begin{array}{c|c|c|c}
\alpha & -\frac 1 2 (\alpha,\alpha) & \mathrm{mult}(\alpha) & \text{main term} \\\hline
(7,7,2) & 10 & 56 & 56.65 \\
(8,10,4) & 20 & 792 & 793.19 \\
(11,11,5) & 30 & 6826 & 6867.52 \\
(11,14,7) & 40 & 44258 & 44975.14
\end{array}
\end{equation}
  (A table of $\mathrm{mult}(\alpha)$ can be found in \cite[p.205]{K}.)
\end{example}

Using the fact
$I_2(x)\sim \frac {e^x}{\sqrt{2\pi x}}$, we can see $p_{\sigma}(n+1)\sim \tfrac {e^{\frac {4\pi\sqrt{n}}{\sqrt{23}}}}{n^{\frac 54}23^{\frac 14}\sqrt{2}}$. It is interesting to compare it with $p(n)\sim \tfrac {e^{\pi\sqrt{\frac {2n}3}}}{4n\sqrt{3}}$ and to see the deviation from Frenkel's conjecture.

This method can be applied to other hyperbolic Kac--Moody algebras and to other modular forms as shown in \cite{KL}. 


\section{Summary of results on root multiplicities}

 There are recursive formulas for root multiplicities by Peterson (\cite{P}) and Kang (\cite{Ka2}) as well as closed form formulas by Berman and Moody (\cite{BM}) and Kang (\cite{Ka2}). 

 A recursive formula for root multiplicities of hyperbolic or Lorentzian Kac--Moody algebras assumes knowledge of the root multiplicities corresponding to subalgebras of finite or affine type.

 The known closed form formulas of Berman and Moody and Kang require a substantial amount of information for their application, such as representations and root multiplicities of subalgebras. Moreover they give answers for root multiplicities one at a time, with no general formulas or effective bounds on multiplicities.

 There are many partial results for root multiplicities of hyperbolic and Lorentzian Kac-Moody algebras by applying the formulas of   Peterson, Berman and Moody, and Kang with additional external data as in Feingold and Frenkel \cite{FF}, Kac, Moody and Wakimoto \cite{KMW},  Benkart, Kang and Misra \cite{BKM1},  Kang and Melville \cite{KM2}, Klima and Misra \cite{KMi}, Hontz and Misra \cite{HM1}, Kleinschmidt \cite{Kl2}, Bauer and Bernard \cite{BB}.
However these  results  do not suggest any unified approach to computing root multiplicities.

The only result which suggests a unified viewpoint on root multiplicities for hyperbolic and Lorentzian Kac--Moody algebras is a conjecture of Frenkel. In this approach, it is still an open problem to formulate  precise bounds for root multiplicities for hyperbolic and other indefinite Kac--Moody algebras. Such a conjecture may emerge from connections of hyperbolic Kac--Moody algebras to automorphic forms, which would make it possible to use powerful, analytic tools.

\noindent {\bf Acknowledgments} We are grateful to I. Frenkel and S.-J. Kang  for their helpful comments and encouragement. We thank A. Kleinschmidt for providing his thesis and other computations as well as useful references. The third named author stayed at the Institute for Computational and Experimental Research in Mathematics (ICERM) while this work was in progress, and he wishes to acknowledge its stimulating research environment.


\bibliographystyle{amsalpha}

\end{document}